\documentclass[12pt]{article}
\usepackage{amsmath}
\usepackage{amssymb}
\usepackage{graphicx}
\usepackage[cp1251]{inputenc}
\usepackage[russian]{babel}
\newcommand{\il}[2]{\int\limits_{#1}^{#2}}
\newcommand{\ilp}[1]{\int\limits_{#1}^{+\infty}}

\newcommand{\ph}{\phantom{a}}
\newcommand{\phh}{\phantom{aaa}}
\newcommand{\ilpp}{\ilp{t_0}}
\newcommand{\sist}[2]{\left\{
\begin{array}{l}
{#1}\\
\ph\\
{#2}
\end{array}
\right.}

\begin{document}

\vskip 20pt

MSC 34C10

\vskip 20pt

\centerline{\bf Oscillatory and non oscillatory criteria for the  systems  of two linear}
 \centerline{\bf  first  order two by two dimensional  matrix ordinary differential equations}

\vskip 20 pt

\centerline{\bf G. A. Grigorian}
\centerline{\it Institute  of Mathematics NAS of Armenia}
\centerline{\it E -mail: mathphys2@instmath.sci.am}
\vskip 20 pt

\noindent
Abstract.  The Riccati equation method is used for study the oscillatory and non oscillatory behavior of solutions of systems of two first order linear two by two dimensional matrix differential equations. An integral and an interval oscillatory criteria are obtained. Two non  oscillatory criteria are obtained as well. On an example one of the obtained oscillatory criteria is compared with some well known  results.
\vskip 20 pt


\vskip 20 pt

\centerline{\bf \S\hskip 3pt 1. Introduction}

\vskip 20 pt

 Let   $P(t)\equiv \bigl(p_{jk}(t)\bigr)_{j,k=1}^2, \ph Q(t) \equiv diag \{q_1(t), q_2(t)\}, \ph R(t)\equiv \bigl(r_{jk}(t)\bigr) _{j,k=1}^2, \linebreak   S(t)\equiv \bigl(s_{jk}(t)\bigr)_{j,k=1}^2$                                             be real valued continuous matrix functions on  $[t_0;+\infty)$. Consider the matrix linear system
$$
\sist{\Phi'= P(t)\Phi + Q(t)\Psi;}{\Psi' = R(t)\Phi + S(t)\Psi, \phh t\ge t_0.} \eqno (1.1)
$$
Here $\Phi=\Phi(t)\equiv \bigl(\phi_{jk}(t)\bigr)_{j,k=1}^2,\ph  \Psi=\Psi(t)\equiv\bigl(\psi_{jk}(t)\bigr)_{j,k=1}^2$ are unknown continuously differentiable matrix  functions on  $[t_0;+\infty)$.

{\bf Remark 1.1.} {\it The general case $Q(t)\equiv S(t) diag\{q_1(t), q_2(t)\} S^{-1}(t), \ph t\ge t_0$, where $S(t)$ is an invertible continuously differentiable on $[t_0;+\infty)$ matrix  function, can be reduced to the case $Q(t) \equiv diag\{q_1(t), q_2(t)\}, \ph t\ge t_0,$ of the system (1.1) by the linear transfor-\linebreak mation
$$
\Phi = S(t) \Phi_1, \phh \Psi = S(t) \Psi_1, \phh t\ge t_0.
$$
in (1.1).
}

{\bf Definition 1.1.} {\it A solution  $(\Phi(t), \Psi(t))$ of the system (1.1) is called oscillatory if   $\det \Phi(t)$  has arbitrary large zeroes, otherwise it is called non oscillatory.}

{\bf Definition 1.2.} {\it A solution  $(\Phi(t), \Psi(t))$ of the system (1.1) is called oscillatory on the interval $[t_1;t_2], \ph (t_0 \le t_1 < t_2 < +\infty)$ if  $\det \Phi(t)$  has at least one zero on $[t_1;t_2]$.}

{\bf Definition 1.3}. {\it A solution  $(\Phi(t), \Psi(t))$ of the system (1.1) is called prepared \linebreak (or preferred) if  $\Phi^*(t)\Psi(t) = \Psi^*(t)\Phi(t), \ph t\ge t_0$, where   $*$  is the transpose sign.}

{\bf Definition  1.4}. {\it The system (1.1) is called oscillatory,
if its all prepared solutions are  oscillatory.}

{\bf Definition  1.5}. {\it The system (1.1) is called oscillatory on the interval $[t_1;t_2], \linebreak (t_0 \le t_1 < t_2 < +\infty)$ if its all prepared solutions are oscillatory on the interval $[t_1;t_2]$.}

Study of questions of oscillation and non oscillation of solutions of linear systems of matrix equations, in particular of the system
(1.1), is an important problem of qualitative theory of differential equations and many works are devoted to them (see for example [1~ - 7]). In most of cases in the works  [1 - 7] and others on the matrix  coefficients of the system are imposed conditions ensuring some symmetry property of corresponding matrix Riccati equation (the hamiltonian systems), namely if  $Y(t)$ is a solution to corresponding Riccati equation then the transposed matrix  function $Y^*(t)$ is a solution of the last one as well.
In this work we study the conditions on the coefficients of the system
 (1.1),  for which the last one has oscillatory and non oscillatory solutions. We impose conditions on the coefficients of the system (1.1) for which the hamiltonian structure of it can  not be kept.

 \vskip 20 pt

\centerline{\bf \S\hskip 3pt 2. Auxiliary propositions}

 \vskip 20 pt

In this paragraph we prove two lemmas  and  represent a lemma and a theorem,  proved in other works. They will be used in the next paragraph for proving oscillatory and non oscillatory criteria for the system (1.1).

In what follows the solutions of equations and systems of equations we will assume  real valued. In the system (1.1) make a change
$$
\Psi = Y(t)\Phi, \phh t\ge t_0, \eqno (2.1)
$$
where $Y(t)$  is a continuously differentiable matrix  function of dimension  $ 2\times 2$ on  $[t_0;+\infty)$.  We will get:
$$
\sist{\Phi'=[P(t) + Q(t) Y(t)] \Phi;}{[Y'(t) + Y(t) Q(t) Y(t) + Y(t) P(t) - S(t) Y(t) - R(t)]\Phi = 0,\phh t\ge t_0.} \eqno (2.2)
$$
Consider the matrix Riccati equation
$$
Y' + Y Q(t) Y + Y P(t) - S(t) Y - R(t) = 0, \phh t\ge t_0, \eqno (2.3)
$$
where  $Y=\bigl(y_{jk}(t)\bigl)_{j,k=1}^2$. From (2.2) is seen that if  $Y_1(t)$  is a solution of Eq.  (2.3) on  $[t_1;t_2) \linebreak (t_0\le t_1 < t_2 \le +\infty)$, then  $(\Phi_1(t), Y_1(t)\Phi_1(t))$ is a solution to the system (1.1) on  $[t_1;t_2)$, where $\Phi_1(t)$ is any solution to matrix equation
$$
\Phi' = [P(t) + Q(t)Y_1(t)]\Phi, \phh t\in [t_1;t_2). \eqno (2.4)
$$
Obviously on the strength of  (2.1) and (2.2) if  $(\Phi(t), \Psi(t))$ is a solution of the system (1.1) and  $\det \Phi(t) \ne~0, \ph t\in [t_1;t_2)$, then $Y(t) \equiv \Psi(t)\Phi^{-1}(t)$ is a solution to Eq. (2.3)  on  $[t_1;t_2)$.
Let  $Y_0(t)$  be a solution to Eq. (2.3) on  $[t_1;t_2)$.

{\bf Definition 2.1.} {\it We will say that $[t_1;t_2)$  is a maximum existence interval for $Y_0(t)$, if   $Y_0(t)$  cannot be continued to the right of $t_2$ as a solution of Eq. (2.3).}

{\bf Lemma 2.1.} {\it Let $Y_0(t)$  be a solution of Eq. (2.3) on  $[t_1;t_2)$,  and let  $t_2< +\infty$.
Then $[t_1;t_2)$ cannot be the maximum existence interval for $Y_0(t)$ provided the function\linebreak $f(t)\equiv \il{t_1}{t} tr [Q(\tau) Y_0(\tau)]d\tau, \ph t\in [t_1;t_2)$, is bounded from below on $[t_1;t_2)$.}

Proof. Let $\Phi_0(t)$  be a solution to the equation
$$
\Phi' = [P(t) + Q(t) Y_0(t)] \Phi, \phh t\ge t_0, \eqno (2.5)
$$
with  $\Phi_0(t_1) \ne 0$. Then by Liouville formula
$$
\det \Phi_0(t) = \det\Phi_0(t_1)\exp\biggl\{\il{t_1}{t} tr \bigl[P(\tau) + Q(\tau) Y_0(\tau)\bigr] d\tau \biggr\} \ne 0, \phh t\in [t_1;t_2). \eqno (2.6)
$$
Recall that for any solution $\Phi_0(t)$ of the linear matrix equation
$$
\Phi' = A(t)\Phi, \phantom{aaa} t\ge t_0,
$$
where $A(t)$ is a square continuous matrix function, the Liuville's theorem states that (the Liuville's formula)
$$
det \Phi_0(t) = det\Phi_0(t_0) \exp\biggl\{\int\limits_{t_0}^t tr(A(\tau)) d \tau\biggl\}
$$
(see [8], p. 47, Theorem 1,2).
Let  $(\widetilde{\Phi}(t), \widetilde{\Psi}(t))$ be the solution of the system (1.1) with  $\widetilde{\Phi}(t_1) = \Phi_0(t_1), \ph \widetilde{\Psi}(t_1) = Y_0(t_1)\Phi_0(t_1)$. Then by (2.2)  - (2.5) and the uniqueness theorem
$$
\widetilde{\Phi}(t) = \Phi_0(t), \phh \widetilde{\Psi}(t) = Y_0(t) \Phi_0(t), \phh t\in [t_1;t_2). \eqno (2.7)
$$
From the conditions of the lemma and from (2.6) it follows that $|det \Psi_0(t)| \ge \varepsilon, \ph t\in [t_1;t_2),$ for some $\varepsilon > 0$. Then since $det \widetilde{\Psi}_0(t)$ is a continuous function from (2.7)
it follows that $\det \widetilde{\Phi}(t)  \ne 0, \ph t\in [t_1;t_3)$,  for some $t_3 > t_2$.  Therefore $\widetilde{Y}_0(t)\equiv \widetilde{\Psi}(t)\widetilde{\Phi}^{-1}(t)$   is a solution to Eq. (2.3) on  $[t_1;t_3)$. By (2.7) we have  $\widetilde{Y}_0(t) = Y_0(t), \ph t\in [t_1;t_2)$. Hence  $[t_1;t_2)$ is not the maximum existence interval for  $Y_0(t)$.  The lemma is proved.

Let $a(t), \ph b(t), \ph c(t), \ph c_1(t)$ be continuously differentiable functions on  $[t_0;+\infty)$.\linebreak  Consider the Riccati equations
$$
y' + a(t) y^2 + b(t) y + c(t) = 0, \phh t\ge t_0; \eqno (2.8)
$$
$$
y' + a(t) y^2 + b(t) y + c_1(t) = 0, \phh t\ge t_0; \eqno (2.9)
$$

{\bf Theorem 2.1}. {\it Let Eq. (2.9) has the solution $y_1(t)$  on  $[t_1;t_2) \ph (t_0\le t_1 < t_2 \le +\infty)$, and let  $a(t) \ge 0, \ph c(t) \le c_1(t), \ph t\in [t_1;t_2)$. Then for each  $y_{(0)} \ge y_1(t_0)$ Eq. (2.8)
has the solution $y_0(t)$  on  $[t_1;t_2)$   with  $y_0(t_0) = y_{(0)}$, and $y_0(t) \ge y_1(t), \ph t\in [t_1;t_2)$.}

A proof for a more general theorem is presented in [9] (see also [10]).

Let us write Eq. (2.3) in the expanded form. We have:
$$
\left\{
\begin{array}{l}
y_{11}' + q_1(t) y_{11}^2 + a_{11}(t) y_{11} + q_2(t) y_{12} y_{21} + p_{21}(t) y_{12} - s_{12}(t) y_{21} - r_{11}(t) =0;\\
y_{22}' + q_2(t) y_{22}^2 + a_{22}(t) y_{22} +q_1(t) y_{12} y_{21} + p_{12}(t) y_{21}- s_{21}(t) y_{12} - r_{22}(t) = 0;\\
y_{12}' + [q_1(t) y_{11} + q_2(t) y_{22} + a_{21}(t)]y_{12} + p_{12}(t) y_{11}- s_{12}(t) y_{22} - r_{12}(t)=0;\\
y_{21}' + [q_1(t) y_{11} + q_2(t) y_{22} + a_{12}(t)]y_{21} + p_{21}(t) y_{22}  - s_{21}(t) y_{11} - r_{21}(t)=0,
\end{array}
\right. \phantom{aa} \eqno (2.10)
$$
where  $a_{jk}(t) \equiv p_{jj}(t) - s_{kk}(t), \ph j, k =1,2, \ph t\ge t_0$.
Denote:\\

\noindent
$
I_k(\tau;t) \equiv \il{\tau}{t}\exp\biggl\{- \il{s}{t} a_{kk}(\zeta)d\zeta\biggr\} r_{kk}(s) d s, \ph t \ge \tau \ge t_0, \ph k=1,2.
$

{\bf Lemma 2.2.} {\it Let the following conditions hold

\noindent
A) $q_k(t) \ge 0, \ph k =1,2, \ph r_{12}(t) \ge 0 \ph (\le 0), \ph r_{21}(t) \le 0 \ph (\ge 0), \ph p_{21}(t) \ge 0 \ph (\le 0), \linebreak s_{12}(t) \ge 0 \ph (\le 0), \ph t\ge t_0$;

\noindent
B) there exist infinitely large sequences $\xi_{0,k} = t_0 < \xi_{1,k} < \dots  < \xi_{m,k} < \dots\linebreak \ph (k=1,2)$ such that
$$
\il{\xi_{m,k}}{t}\exp\biggl\{\il{\xi_{m,k}}{\tau}\bigr[a_{kk}(s) + q_k(s) I_k(\xi_{m,k};s)\bigr]d s\biggr\} r_{kk}(\tau)d\tau \ge 0, \ph t\in [\xi_{m,k};\xi_{m+1,k}), \ph m=0,1,2, \dots,
$$
$k=1,2.$
Then for each  $y_{kk,0} > 0, \ph k=1,2, \ph y_{12,0} \le 0 \ph (\ge 0) \ph  y_{21,0} \ge 0 \ph (\le 0)$ \linebreak Eq. (2.3) has the solution  $Y_0(t) \equiv \bigl(y_{jk}^0(t)\bigl)_{j,k=1}^2$  on  $[t_0;+\infty)$,  satisfying the initial conditions $y_{jk}^0(t_0) = y_{jk,0}, \ph j,k = 1,2,$
and
$$
\det Y_0(t) > 0, \phh t\ge t_0. \eqno (2.11)
$$
}

Proof. Show that
$$
y_{kk}^0(t) > 0, \phh t\in [t_0;T), \phh k=1,2, \eqno (2.12)
$$
where $[t_0;T)$ is the maximum existence interval for $Y_0(t)$.  Suppose that it is not so. Then from the initial conditions is seen that
$$
y_{kk}^0(t) > 0, \phh t\in [t_0;T_1), \eqno (2.13)
$$
$$
y_{11}^0(T_1) y_{22}^0(T_1) = 0, \eqno (2.14)
$$
for some  $T_1 \in (t_0;T)$. By virtue of the third and fourth equations of the system (2.10) we have:
$$
y_{12}^0(t) = \exp\biggl\{- \il{t_0}{t}[q_1(\tau) y_{11}^0(\tau) + q_2(\tau) y_{22}^0(\tau) + a_{21}(\tau)]d\tau\biggr\}\biggl[y_{12}^0(t_0) - \phantom{aaaaaaaaaaaaaaaaaaaaa}
$$
$$
- \il{t_0}{t}\exp\biggl\{\il{t_0}{\tau}\bigl[q_1(s)y_{11}^0(s) + q_2(s) y_{22}^0(s) + a_{21}(s)\bigr] d s\biggr\}\times \phantom{aaaaaaaaaa}
$$
$$
\phantom{aaaaaaaaaaaaaaaaaaaaa}\times \biggl(p_{12}(\tau) y_{11}^0(\tau) - s_{12}(\tau) y_{22}(\tau) - r_{12}(\tau)\biggr) d\tau \biggr], \ph t\in [t_0;T); \eqno (2.15)
$$

$$
y_{21}^0(t) = \exp\biggl\{- \il{t_0}{t}[q_1(\tau) y_{11}^0(\tau) + q_2(\tau) y_{22}^0(\tau) + a_{12}(\tau)]d\tau\biggr\}\biggl[y_{21}^0(t_0) - \phantom{aaaaaaaaaaaaaaaaaaaaa}
$$
$$
- \il{t_0}{t}\exp\biggl\{\il{t_0}{\tau}\bigl[q_1(s)y_{11}^0(s) + q_2(s) y_{22}^0(s) + a_{12}(s)\bigr] d s\biggr\}\times \phantom{aaaaaaaaaa}
$$
$$
\phantom{aaaaaaaaaaaaaaaaaaaaa}\times \biggl(p_{21}(\tau) y_{22}^0(\tau) - s_{21}(\tau) y_{11}(\tau) - r_{21}(\tau)\biggr) d\tau \biggr], \ph t\in [t_0;T); \eqno (2.16)
$$
From here from the conditions of lemma and from (2.13) it follows that
$$
y_{12}^0(t) \ge 0 \ph (\le 0), \phh y_{21}^0(t) \le 0 \ph (\ge0), \phh t\in [t_0;T_1). \eqno (2.17)
$$
Consider the Riccati equations
$$
y' + q_k(t) y^2 + a_{kk}(t) y - r_{kk}(t) = 0, \phh t\ge t_0, \eqno (2.18_k)
$$
$$
y' + q_k(t) y^2 + a_{kk}(t) y +  \mathcal{L}_k(y_{k,3-k}^0(t),y_{3-k,k}^0(t),t) = 0, \phh t\ge t_0, \eqno (2.19_k)
$$
$k=1,2$,
where $\mathcal{L}_k(u,v,t) \equiv q_{3-k}(t)u v + p_{3-k,k}(t) u - s_{k,3-k}(t)v - r_{kk}(t), \ph u,v \in R,\linebreak t\ge t_0, \ph k=1,2.$
From the conditions A) of lemma and from (2.17) it follows that
$$
\mathcal{L}_k(y_{k,3-k}^0(t),y_{3-k,k}^0(t),t)  \le - r_{kk}(t), \phh t\in [t_0;T_1), \phh k=1,2. \eqno (2.20)
$$
Let  $y_k(t)$ be the solution of Eq. $(2.18_k)$ with   $y_k(t_0) = y_{kk}^0(t_0) > 0, \ph (k=1,2)$.               Then on the strength of Theorem 4.1 of work [10] from the conditions  B) of lemma it follows that $y_k(t)$  exists on $[t_0;T)$  and
$$
y_k(t) > 0, \phh t\in [t_0;T),\phh k=1,2. \eqno (2.21)
$$
Obviously by (2.10) the function $y_{kk}^0(t)$ is a solution to Eq. $(2.19_k)$  on  $[t_0;T), \linebreak (k=1,2)$.
Then by virtue of Theorem 2.1 from (2.20) and (2.21) it follows that  \linebreak  $y_{kk}^0(t) \ge y_k(t) > 0, \ph t\in [t_0;T_1], \ph k=1,2$, which contradicts (2.14). The obtained contradiction proves (2.12). Show that $T=+\infty$. From the conditions $q_k(t) \ge 0,  \linebreak  t\ge t_0, \ph  k=1,2$  (a part of A)), and from  (2.12) it follows that
$$
\il{t_0}{t} tr [Q(\tau) Y_0(\tau)] d\tau \ge 0, \phh t\in [t_0;T). \eqno (2.22)
$$
Suppose  $T< +\infty$. Then  by Lemma 2.1 from (2.22) it follows that  $[t_0;T)$ is not the maximum existence interval for $y_0(t)$.   The obtained contradiction shows that $T=+\infty$. From here, from the conditions   A) of lemma, from (2.12), (2.15) and (2.16) it follows (2.11). The lemma is proved.

{\bf Remark  2.1.} {\it The conditions B)  of Lemma 2.2 are  satisfied if in particular \linebreak $r_{kk}(t) \ge 0, \ph t\ge t_0, \ph k=1,2.$}

{\bf Lemma 2.3.} {\it Let Eq. (2.8) has a solution on $[t_1;+\infty)$  for some $t_1 \ge t_0$ , and let  $a(t) \ge0, \ph c(t) \ge 0, \ph t\ge t_0, \ph \ilpp a(\tau)\exp\biggl\{- \il{t_0}{\tau}b(s) d s\biggr\}d\tau = +\infty.$
Then Eq. (2.8) has a positive solution on  $[t_1;+\infty)$.}

The proof is presented in [11].

\vskip 20pt

\centerline{\bf \S \hskip 3pt 3. Oscillatory and non oscillatory criteria}

\vskip 20pt

Denote:
$$
F_k(t) \equiv \sist{r_{kk}(t) - (p_{3-k,k}(t) - s_{k,3-k}(t))^2/(4q_{3-k}(t)), \ph q_{3-k}(t) \ne 0;}{r_{kk}(t), \phantom{aaaaaaaaaaaaaaaaaaaaaaaaaaaaaa} q_{3-k}(t) = 0,}
$$
$t\ge t_0, \ph k=1,2$. Let $j(\in \{1,2\})$  be fixed.
Consider the Riccati equation
$$
y' + q_j(t)y^2 + a_{jj}(t) y - F_j(t) = 0, \phh t\ge t_0. \eqno (3.1)
$$
The solutions $y(t)$ of this equation  existing on some interval  $[t_1;t_2) \linebreak (t_0 \le t_1 < t_2 \le +\infty)$,  are connected with the solutions  $(\phi(t), \psi(t))$  of the system of scalar equations
$$
\sist{\phi' = p_{jj}(t) \phi + q_j(t) \psi;}{\psi' = F_j(t) \phi + s_{jj}(t) \psi, \ph t\ge t_0,} \eqno (3.2)
$$
by relations  (see [12])
$$
\phi(t) = \phi(t_1)\exp\biggl\{\il{t_1}{t}\bigl[q_j(\tau) y(\tau) + p_{jj}(\tau)\bigr]d\tau\biggr\}, \ph \phi(t_1) \ne 0, \ph \psi(t) = y(t) \phi(t), \eqno (3.3)
$$
$t\in [t_1;t_2).$

{\bf Definition 3.1.} {\it  The system  $(3.2)$ is called oscillatory if for its each solution  $(\phi(t), \psi(t))$ the function  $\phi(t)$  has arbitrary large zeroes.}

{\bf Definition 3.2.} {\it  The system  $(3.2)$ is called oscillatory on the interval $[t_1;t_2]$ if for its each solution  $(\phi(t), \psi(t))$ the function  $\phi(t)$  has at least one zero on $[t_1;t_2]$.}

{\bf Theorem 3.1}. {\it  Let the following conditions be satisfied:

\noindent
I)  $q_k(t) \ge 0, \ph k=1,2,$ and if $q_{3-j}(t) =0$, then $p_{3-j,j}(t) = s_{j,3-j}(t), \ph t\ge t_0;$

\noindent
II) the system (3.2) is oscillatory.

\noindent
Then the system (1.1) is oscillatory.}

Proof.  Let $(\Phi(t),\Psi(t))$ be a prepared solution to the system (1.1). Suppose that $(\Phi(t),\Psi(t))$  is not oscillatory.   Then $\det \Phi(t) \ne 0, \ph t\ge T$,  for some $T \ge t_0$.
Let \linebreak $Y_0(t) \equiv \bigl(y_{jk}^0(t)\bigr)_{j,k = 1}^2 = \Psi(t) \Phi^{-1}(t), \ph t\ge T$.
By (2.1)  $Y_0(t)$ is a solution of Eq. (2.3) on  $[T;+\infty)$. Then by (2.10) $y_{jj}^0(t)$ satisfies to the following Riccati equation
$$
y' + q_j(t) y^2 + a_{jj}(t) y + \mathcal{L}_j(y_{j,3-j}^0(t),y_{3-j,j}^0(t),t)= 0, \phh t\ge T \eqno (3.4)
$$
(the definition of $\mathcal{L}_j$\ph see below $(2.19_k)$).  Since $(\Phi(t),\Psi(t))$   is a prepared solution we have  $Y_0(t) = Y_0^*(t), \ph t\ge T$.
From here and from the conditions I)  of theorem it follows that
$$
\mathcal{L}_j(y_{j,3-j}^0(t),y_{3-j,j}^0(t),t) \ge F_j(t), \phh t\ge T. \eqno (3.5)
$$
Consider the Riccati equation
$$
y' + q_j(t) y^2 + a_{jj}(t) y  - F_j(t) = 0, \phh t\ge T. \eqno (3.6)
$$
Let  $y_j(t)$  be its solution with  $y_j(t) \ge y_{jj}^0(T)$.
Then using Theorem~2.1 by applying (3.5) to the equations (3.4) and (3.6) we will conclude that $y_j(t)$    exists on  $[T;+\infty)$.   Therefore by (3.1) - (3.3) the functions
$$
\phi_j(t) =\exp\biggl\{\il{T}{t}\bigl[q_j(\tau) y(\tau) + p_{jj}(\tau)\bigr]d\tau\biggr\}, \ph  \psi_j(t) = y_j(t) \phi_j(t), \ph t\ge T
$$
form the solution $(\phi_j(t), \psi_j(t))$ of the system (3.2) on $[T;+\infty)$, which can be continued on $[t_0;+\infty)$  as a solution of the system (3.2). It is evident that $\phi_j(t)$  has no arbitrary large zeroes which contradicts  II). The theorem is proved.

By analogy can be proved

{\bf Theorem 3.2}. {\it  Let the following conditions be satisfied:

\noindent
I$^*$)  $q_k(t) \ge 0, \ph k=1,2,$ and if $q_{3-j}(t) =0$, then $p_{3-j,j}(t) = s_{j,3-j}(t),  \ph t\in [t_1;t_2] \\ (t_0 \le t_1 < t_2 < +\infty);$

\noindent
II$^*$) the system (3.2) is oscillatory on the interval $[t_1;t_2]$.

\noindent
Then the system (1.1) is oscillatory on the interval $[t_1;t_2]$.}

{\bf Remark 3.1}. {\it The restrictions I) on $Q(t)$ in Theorem 3.1 means that  $Q(t)$ is nonne-\linebreak gative definite meanwhile  in the works  [1 - 7] and others the corresponding coefficient is positive definite.}

{\bf Remark 3.2}. {\it Suppose  $p_{12}(t)= - s_{21}(t), \ph p_{12}(t) = - s_{21}(t) , \ph a_{12}(t) = a_{21}(t), \linebreak r_{12}(t) = r_{21}(t), \ph t\ge t_0$.  Then by  (2.10), if  $Y_0(t)$  is a solution of Eq. (2.3) on\linebreak some interval $[t_0;t_1)$, then $Y_0^*(t)$  is a solution of Eq. (2.3) on $[t_0;t_1)$ too. On the strength\linebreak of the uniqueness theorem from here it follows that if  $Y_0(t_0)=Y_0^*(t_0)$, then\linebreak   $Y_0(t) = Y_0^*(t), \ph t\in [t_0;t_1)$. Therefore taking into account  (2.1) we conclude that if $(\Phi(t), \Psi(t))$  is a solution of the system (1.1)  with
$\det \Phi(t_0) \ne 0, \ph \Phi^*(t_0)\Psi(t_0)= \Psi^*(t_0) \Phi(t_0)$,  then  $\Phi^*(t)\Psi(t)= \Psi^*(t) \Phi(t), \ph t\in [t_0;t_1).$ Obviously the last equality will be satisfied on the whole interval $[t_0;+\infty)$, provided we additionally require that  $P(t), \ph Q(t), \ph R(t)$ and $S(t)$  be analytical functions on the some domain  of complex plane containing the half line  $[t_0;+\infty)$.
From the given restrictions  above on  $P(t), \ph Q(t), \ph R(t)$ and  $S(t)$ is seen that the system (1.1) can be not hamiltonian.  So the system (1.1) can have prepared solution not only in the case when it is hamiltonian but also in the other cases}.

Example 3.1. Consider the matrix equation
$$
\Phi'' + K(t) \Phi = 0, \phantom{aaa} t\ge t_0. \eqno (3.7)
$$

\noindent
where
$
 K(t) \equiv\left( \begin{array}{c} a_1 \sin \mu_1 t + a_2 \sin \mu_2 t \phantom{aaa} \frac{b \cos \mu t}{t^\alpha}\\
\frac{b \cos \mu t}{t^\alpha}\phantom{aaa} a_1 \sin \mu_1 t + a_2 \sin \mu_2 t \end{array}\right),\ph a_1, \ph a_2, \ph \alpha \ph \mu, \ph \mu_1, \ph \mu_2
$
are some real nonzero constants and  $\alpha > 1, \ph \mu_1/ \mu_2$ is irrational. This equation  is equivalent to the system (1.1) with  $P(t)= S(t) \equiv 0, \ph  R(t) \equiv K(t) \ph Q(t) \equiv I$              where $I$  is the identity matrix of dimension $2\times 2$.  Therefore for this equation the system (3.2) has the form
$$
\sist{\phi' = \phantom{aaaaaaaaaaaaaaaaaaaaaaaaaa} \psi;}{\psi' = - (a_1 \sin \mu_1 t + a_2 \sin \mu_2 t) \phi, \ph t\ge t_0.}
$$
which is equivalent to the scalar equation
$$
\phi'' + (a_1 \sin \mu_1 t + a_2 \sin \mu_2 t) \phi = 0, \phh t\ge t_0.
$$
This equation is oscillatory (see  [13], Corollary 3.4). Therefore the last system is oscillatory too.
By virtue of Theorem 3.1 from here it follows that Eq. (3.7) is oscillatory.
The eigenvalues $\lambda_\pm(t)$   of the matrix  $K(t)$   are equal
$$
\lambda_\pm(t) =  a_1 \sin \mu_1 t + a_2 \sin \mu_2 t \hskip 2pt \pm \hskip 2pt   \frac{|b \cos \mu t|}{t^\alpha}, \phh t\ge t_0.
$$
From here is seen that the Theorems 5,  6  of work  [1], and the Theorems 1, 2 and 3 of work [14]
are not applicable to  Eq. (3.7).
The remaining theorems of these works  and the results of works  [2 -7] are not explicit for applying them to Eq. (3.7) (it is hard to guess can we apply them  to Eq. (3.7)).

{\bf Corollary 3.1.} {\it Let the conditions I)  of Theorem 3.1 be satisfied and let

\noindent
III)  $ \ilpp q_j(\tau)\exp\biggl\{-\il{t_0}{\tau} a_{jj}(s)ds\biggr\}d\tau = \ilpp[-F_{j}(\tau)]\exp\biggl\{\il{t_0}{\tau}a_{jj}(s) d s\biggr\}d \tau =+\infty$.

\noindent
Then the system (1.1) is oscillatory.}

Proof. On the strength of Theorem 3.1 it is enough to show that the system (3.2) is oscillatory. Suppose that the system  (3.2) is not oscillatory. Then by (3.1) - (3.3) Eq. (3.1) has a solution on  $[t_1;+\infty)$ for some  $t_1\ge t_0$. Set  $W(t) \equiv  -F_{j}(t) \exp\biggl\{\il{t_1}{t} a_{jj}(\tau) d\tau\biggr\}, \ph t\ge t_1.$ In Eq. (3.1) make the change
$$
y= z \exp\biggl\{- 2\il{t_1}{t} a_{jj}(\tau) d\tau\biggr\}, \phh t\ge t_1.
$$
We will come to the equation
$$
z' + U(t) z^2 + W(t) = 0, \phh t\ge t_1, \eqno (3.8)
$$
where  $U(t) \equiv q_j(t) \exp \biggl\{-\il{t_1}{t} a_{jj}(\tau) d\tau\biggr\}.$
Show that
$$
\ilp{t_1} U(\tau)\exp\biggl\{\il{t_1}{t} U(s) d s \il{t_1}{s} W(\zeta) d \zeta\biggr\} d\tau = +\infty. \eqno (3.9)
$$
On the strength of III) we have: $\il{t_1}{t}W(\tau) d \tau = - \il{t_1}{t} F_{j}(\tau) \exp \biggl\{\il{t_1}{t} a_{jj}(s) d s\biggr\} d \tau \ge 0, \ph t\ge t_2$, for some  $t_2 \ge t_1$. By  III)
from here it follows (3.9). In Eq. (3.8) make the change
$$
z = u - \il{t_1}{t} W(\tau) d\tau, \phh t\ge t_1.
$$
We will get
$$
u' + U(t) u^2 - 2U(t)\il{t_1}{t} W(\tau) d\tau u + U(t)\biggl[ \il{t_1}{t} W(\tau) d\tau\biggr]^2 = 0, \phh t\ge t_1. \eqno (3.10)
$$
Since by assumption Eq. $(3.1_j)$ has a solution on $[t_1;+\infty)$,  from the above substitutions is seen that Eq. (3.10) has a solution on $[t_1;+\infty)$ as well. By virtue of Lemma 2.3 from here from (3.9) and from the inequalities  $q_j(t) \ge 0, \ph U(t)\biggl[ \il{t_1}{t} W(\tau) d\tau\biggr]^2 \ge 0, \ph t\ge t_1$
it follows that Eq. (3.10) has a positive solution $u_0(t)$ on   $[t_1;+\infty)$. Then  $z_0(t) \equiv u_0(t) - \il{t_1}{t} W(\tau) d\tau$  is a solution to Eq. (3.8) such that
$$
z_0(t) > \il{t_1}{t} W(\tau) d\tau, \phh t\ge t_1. \eqno (3.11)
$$
From  (3.8) it follows that
$$
z_0(t) = z_0(t_1) - \il{t_1}{t} U(\tau) z_0^2(\tau) d \tau - \il{t_1}{t} W(\tau) d\tau, \phh t\ge t_1. \eqno (3.12)
$$
From here and from (3.11) we have:
$$
0 \le \il{t_1}{t}U(\tau) z_0^2(\tau) d \tau < z_0(t_1), \phh t\ge t_1. \eqno (3.13)
$$
$(z_0(t_1) = u_0(t_1) > 0)$.
Taking into account  III) from here we will get:\linebreak $\biggl[z_0(t_1) - \il{t_1}{t}U(\tau) z_0^2(\tau) d \tau - \il{t_1}{t} W(\tau) d\tau\biggr]^2 \ge 1, \ph t \ge T,$
for some  $T\ge t_0$. From here and from (3.12) it follows that  $z_0^2(t) \ge 1, \ph t\ge T.$ Therefore by  III) $\ilp{T}U(\tau) z_0^2(\tau) d\tau  \ge \ilp{T} U(\tau) d\tau =\linebreak = + \infty$, which contradicts (3.13). The corollary is proved.

{\bf  Corollary 3.2} {\it Let the conditions $I^*)$ of Theorem 3.2 be satisfied and let

\noindent
IV) $\il{t_1}{t_2} \min \biggl[q_j(\tau)\exp\biggl\{-\il{t_1}{\tau} a_{jj}(s) d s\biggr\}, - F_{j}(\tau) \exp\biggl\{\il{t_1}{\tau} a_{jj}(s) d s\biggr\}\biggr] d \tau \ge \pi$.

\noindent
Then the system (1.1) is oscillatory on the interval $[t_1;t_2]$.
}

Proof. On the strength of Theorem 3.2 it is enough to show that the system  (3.2) is oscillatory on the interval $[t_1;t_2]$. In (3.2) make the changes
$$
\sist{\phi= \exp\biggl\{\int\limits_{t_1}^t p_{jj}(\tau) d \tau\biggr\} \rho \sin \theta;}{\psi= \exp\biggl\{\int\limits_{t_1}^t s_{jj}(\tau) d \tau\biggr\} \rho \cos \theta, \ph t\ge t_0.} \eqno (3.14)
$$
We will get:
$$
\sist{\rho' \sin \theta + \theta'\rho \cos \theta = Q_{j}(t) \rho \cos \theta;}{\rho' \cos \theta - \theta'\rho \sin \theta = R_{j}(t) \rho \sin \theta, \ph t\ge t_0,}  \eqno (3.15)
$$
where $Q_j(t)\equiv \exp\biggl\{-\il{t_1}{t} a_{jj}(\tau) d\tau\biggr\}q_j(t), \ph R_j(t)\equiv  \exp\biggl\{\il{t_1}{t} a_{jj}(\tau) d\tau\biggr\}F_{j}(t), \ph t\ge t_0$ (the function $a_{jj}(t)$ is defined below (2.10)).
This system is equivalent to the system (3.2) in the sense that to each nontrivial solution $(\phi(t),\psi(t))$ of the system (3.2) corresponds the solution $(\rho(t), \theta(t))$ of the system (3.15) with  $\rho(t) >~0, \phantom{a} t\ge t_0$ defined by (3.14).
Let us multiply the first equation of the system (3.15) on $\cos \theta$ and the second one multiply on $\sin \theta$ and subtract from the first obtained equation the second one. We will get:
$$
\theta' \rho = \rho [Q_j(t) \cos ^2 \theta - R_j(t) \sin ^2 \theta], \phh t \ge t_0. \eqno (3.16)
$$
Let $(\phi_0(t), \psi_0(t))$ be a nontrivial solution of the system (3.2) and let $(\rho_0(t), \theta_0(t))$ be the corresponding (by (3.14)) to $(\phi_0(t), \psi_0(t))$ solution of the system (3.15). Then\linebreak $\rho_0(t) \ne 0, \ph t\ge t_0$, and therefore by (3.16) the following equality takes place
$$
\theta_0'(t) = Q_j(t) \cos^2 \theta_0(t) - R_j(t)\sin^2 \theta_0(t) =  \frac{1}{2}\Bigl[Q_j(t) - R_j(t) + (Q_j(t) + R_j(t)) \cos 2 \theta_0(t)\Bigr],
$$
$t\ge t_0$. From here it follows
$$
\theta_0'(t) \ge  \frac{1}{2}\Bigl[Q_j(t) - R_j(t) - |Q_j(t) + R_j(t)|\Bigr] = \min\{Q_j(t), - R_j(t)\}, \phh t\ge t_0.
$$
Let us integrate this inequality from $t_1$ to $t_2$ Taking into account the conditions of the corollary we will get:
$$
\theta_0(t_2) - \theta_0(t_1) \ge  \il{t_1}{t_2} \min \{Q_j(\tau), - R_j(\tau)\} d\tau \ge \pi.
$$
Due to (3.14) from here it follows that $\phi_0(t)$ has at least one zero on $[t_1;t_2]$. The corollary is proved.

{\bf Remark 3.3}. {\it Let $t_0 \le \eta_1 < \zeta_1 < \dots \eta_m < \zeta_m \dots $ be a infinitely large sequence and let the following conditions be satisfied:

\noindent
IV$_m$) $q_k(t) \ge 0,$  and if $q_{3-j}(t) =0$, then $p_{3-j}(t) = s_{j,3-j}(t),  \ph t\in [\eta_m;\zeta_m] , \ph k=1,2;$

\noindent
 $\il{\eta_m}{\zeta_m} \min \biggl[q_j(\tau)\exp\biggl\{-\il{\eta_m}{\tau} a_{jj}(s) d s\biggr\}, - F_{j}(\tau) \exp\biggl\{\il{\eta_m}{\tau} a_{jj}(s) d s\biggr\}\biggr] d \tau \ge \pi, \ph m=1,2,\dots$.

\noindent
Then on the strength of Corollary 3.2 the system (1.1) is oscillatory. From the conditions IV$_m) \ph m=1,2, \dots$ is seen that outside of the set $\bigcup\limits_{m=1}^{+\infty} [\eta_m;\zeta_m]$ the functions $q_1(t)$ and $q_2(t)$ can take values of arbitrary sign and therefore the nonnegative definiteness of $Q(t)$ on $[t_0;+\infty)$
can be broken.
}

{\bf Remark 3.4}. {\it Let $P(t)= S(t)\equiv0, \ph Q(t)= - R(t) \equiv I, \ph t\ge 0$, where $I$ is the identity matrix of dimension $2\times 2$. It is evident that in this particular case the conditions I$^*$) of Corollary 3.2 are satisfied on the arbitrary interval $[t_1;t_2] (\subset [0;+\infty)$) and the condition IV) is fulfilled only if $t_2 - t_1 \ge \pi$. It also is evident that for this case  $(\Phi_0(t), \Psi_0(t))$, where $\Phi_0(t) \equiv diag \{\sin t, \sin t\}, \ph \Psi_0(t) \equiv diag \{\cos t, \cos t\}$, is a prepared solution to the system (1.1). This solution is not oscillatory on $[\varepsilon; \pi - \varepsilon]$ for each $\varepsilon \in (0;\pi)$. Therefore in the inequality IV) we may not replace $\pi$ by a number less than $\pi$.
}

Example 3.2. Consider the system
$$
\sist{\Phi' = \phh K_1(t) \Psi;}{\Psi' = V_1(t) \Phi, \phh t\ge t_0,} \eqno (3.17)
$$
where

$K_1(t) \equiv diag \{ \max\{\sin t, 0\}, \max\{\sin t, 0\}\}, \ph V_1(t) \equiv diag \{ \min\{\sin t, 0\}, \min\{\sin t, 0\}\},$

\noindent
$t\ge t_0$.
Obviously for this system the conditions I$^*$) of Corollary 3.2 are not fulfilled for all $[t_1;t_2] (\subset [t_0;+\infty))$. Therefore Corollary 3.2 cannot be used to establish oscillatory behavior of the system (3.17).
It is easy to verify that for the system (3.17) the conditions of Corollary~3.1 are fulfilled. Therefore the system (3.17)  is oscillatory.

Example 3.3. Consider the system
$$
\sist{\Phi' = \phh K_2(t) \Psi;}{\Psi' = -K_2(t) \Phi, \phh t\ge 0,} \eqno (3.18)
$$
where $K_2(t) \equiv diag\{\lambda \sin t, \lambda \sin t\}, \ph t \ge 0, \ph \lambda \ge \frac{\pi}{2}.$
Obviously the conditions I) of\linebreak Corollary~3.1 for this system are not fulfilled. Therefore it cannot be applied to the system (3.18). It is not difficult to verify that for $t_1 = 2\pi m, \ph t_2= \pi (2 m + 1)$ the conditions of Corollary 3.2 are fulfilled for all $m=1,2, \dots$. Taking into account Remark 3.3 from here we conclude that the system (3.18) is oscillatory.

{\bf Theorem 3.3}. {\it  Let the conditions of Lemma 2.2 be satisfied. Then for each solution  $(\Phi(t), \Psi(t)) \equiv \Bigl(\bigl(\phi_{jk}(t)\bigr)_{j,k =1}^2, \bigl(\psi_{jk}(t)\bigr)_{j,k =1}^2\Bigr)$
of the system (1.1) with  $\det \Phi(t_0) \ne 0, \linebreak y_{11}^0 \equiv \frac{\psi_{11}(t_0) \phi_{22}(t_0) - \psi_{12}(t_0)\phi_{21}(t_0)}{\det \Phi(t_0)} >~0, \ph y_{22}^0 \equiv \frac{\psi_{22}(t_0) \phi_{11}(t_0) - \psi_{21}(t_0)\phi_{12}(t_0)}{\det \Phi(t_0)} > 0, \linebreak y_{12}^0 \equiv \frac{\psi_{12}(t_0) \phi_{11}(t_0) - \psi_{11}(t_0)\phi_{12}(t_0)}{\det \Phi(t_0)} \ge 0 \ph (\le 0), \ph y_{21}^0 \equiv \frac{\psi_{21}(t_0) \phi_{22}(t_0) - \psi_{22}(t_0)\phi_{21}(t_0)}{\det \Phi(t_0)} \le 0 \ph (\ge 0),$
 the equality
$$
sign [\det \Phi(t)] = sign [\det \psi(t)] \ne 0, \phh t\ge t_0. \eqno (3.19)
$$
takes place. Therefore $(\Phi(t),\Psi(t))$ is non oscillatory.}

Proof. On the strength of  Lemma 2.2 Eq. (2.3) has the solution \linebreak
$Y_0(t) \equiv \bigl(y_{jk}(t)\bigl)_{j,k=1,}^2$    on  $[t_0;+\infty)$ with  $y_{jk}(t_0) = y_{jk}^0, \ph j,k =1,2,$ and
$$
\det Y_0(t) > 0, \phh t\ge t_0. \eqno (3.20)
$$
Since by  (2.4) $\Phi(t)$  is a solution to the matrix equation
$$
\Phi' = [P(t) + Q(t) y_0(t)]\Psi, \phh t\ge t_0,
$$
according to  Liouville   formula we have
$$
\det \Phi(t) = \det \Phi(t_0) \exp\biggl\{\il{t_0}{t} tr [P(\tau) + Q(\tau) Y_0(\tau)]d\tau\biggr\} \ne 0, \phh t \ge t_0. \eqno (3.21)
$$
By (2.1) the equality  $\Psi(t) = Y_0(t) \Phi(t), \ph t\ge t_0,$ holds. From here from (3.20) and (3.21) it follows  (3.19). The theorem is proved.

Denote:
$$
\widetilde{I}_1(\tau;t)\equiv\il{\tau}{t}\exp\biggl\{ - \il{s}{t}a_{jj}(\zeta) d\zeta\biggr\} F_j(s) d s,\phantom{aaaaaaaaaaaaaaaaaaaaaaaaaaaaaaaaaaaaaaaa}
$$
$$
\phantom{aaaaaaaaaaaaaaaaaaa}\widetilde{I}_2(\tau;t)\equiv - \il{\tau}{t}\exp\biggl\{ -\il{s}{t}a_{3-j,3-j}(\zeta) d\zeta\biggr\} F_{3-j}(s) d s, \phh t \ge \tau \ge t_0.
$$

{\bf Theorem 3.4}. {\it Let the following conditions be satisfied:

\noindent
C) $q_j(t) \ge 0, \ph q_{3-j}(t) \le 0$ and if $q_j(t)=0$ then $p_{j,3-j}(t) = s_{3-j,j}(t)$, if $q_{3-j}(t) =0$ then $p_{3-j,j}(t) = s_{j,3-j}(t) \ph t\ge t_0;$

\noindent
D) there exists infinitely large sequences $\xi_{0,k} = t_0 < \xi_{1,k} < \dots < \xi_{m,k} < \dots,   \ph k=1,2.$ such that

\noindent
D$_1$) $ \il{\xi_{m,1}}{t}\exp\biggl\{ \il{\xi_{m,1}}{\tau}\Bigl[a_{jj}(s) + q_j(s)\widetilde{I}_1(\xi_{m,1};s)\Bigr]d s\biggr\} F_j(\tau) d\tau \ge 0, \phh t\in [\xi_{m,1}; \xi_{m+1,1}),$

\noindent
D$_2$) $ \il{\xi_{m,2}}{t}\exp\biggl\{- \il{\xi_{m,2}}{\tau}\Bigl[a_{3-j,3-j}(s) + q_{3-j}(s)\widetilde{I}_2(\xi_{m,2};s)\Bigr]d s\biggr\} F_{3-j}(\tau) d\tau \le 0, \phh t\in [\xi_{m,2}; \xi_{m+1,2}),$

\noindent
$m=1,2,\dots.$

\noindent
Then for each prepared solution $(\Phi(t), \Psi(t)) \equiv \Bigl( \bigl(\phi_{jk}(t)\bigr)_{j,k =1}^2, \bigl(\psi_{jk}(t)\bigr)_{j,k =1}^2\Bigr)$ of the system (1.1) with \ph $det \Phi(t_0) \ne 0, \ph y_{11}^0 \equiv\frac{\psi_{11}(t_0) \phi_{22}(t_0) - \psi_{12}(t_0)\phi_{21}(t_0)}{\det \Phi(t_0)}\ge 0, \\ \phantom{aaaaaaaaaaaaaaaaaaaaaaaaaaaaaaaaaaa}y_{22}^0 \equiv \frac{\psi_{22}(t_0) \phi_{11}(t_0) - \psi_{21}(t_0)\phi_{12}(t_0)}{\det \Phi(t_0)} \le 0$ the inequality
$$
det \Phi(t) \ne 0, \phh t\ge t_0, \eqno (3.22)
$$
takes place. Therefore $(\Phi(t), \Psi(t))$ is non oscillatory. Moreover if $y_{11}^0 > 0, \ph y_{22}^0 < 0$, then
$$
sign \hskip 2pt det \Phi(t) = - sign \hskip 2pt det \Psi(t) \ne 0, \phh t\ge t_0. \eqno (3.23)
$$
}

Proof. Let $Y(t) \equiv \bigl(y_{jk}(t)\bigl)_{j,k=1}^2$ be the solution of Eq. (2.3) with $Y(t_0) = \Psi(t_0) \Phi^{-1}(t_0)$, where $(\Phi(t), \Psi(t))$ is a prepared solution to the system (1.1), satisfying the conditions of the theorem, and let $[t_0;T)$ be the maximum existence interval for $Y(t)$. Show that
$$
T = +\infty. \eqno (3.24)
$$
By (2.10) $y_{jj}(t)$ and $y_{3-j,3-j}(t)$ are solutions to the equations
$$
y' + q_j(t) y^2 + a_{jj}(t) y +  \mathcal{L}_j(y_{j,3-j}(t), y_{3-j,j}(t),t) = 0, \ph t\in[t_0;T), \eqno (3.25)
$$
$$
y' - q_{3-j}(t) y^2 + a_{3-j,3-j}(t) y +  \mathcal{L}_{3-j}(y_{3-j,j}(t), y_{j,3-j}(t),t) = 0, \ph t\in[t_0;T), \eqno (3.26)
$$
respectively. From the conditions C) it follows that the following inequalities are satisfied:
$$
\mathcal{L}_j(X,X,t) \le F_j(t), \phh \mathcal{L}_{3-j}(X,X,t) \le - F_{3-j}(t), \phh X\in R, \ph t\ge t_0. \eqno (3.27)
$$
(for $q_j(t) \ne 0 \ph (q_{3-j}(t) \ne 0)$ the $F_j(t) \ph (F_{3-j}(t))$ is the maximum for the quadratic trinomial $\mathcal{L}_j(X,X,t) \ph (\mathcal{L}_{3-j}(X,X,t))$ of variable $X\in R$). Show that
$$
det \Phi(t) \ne 0, \phh t\in [t_0;T). \eqno (3.28)
$$
By (2.4) $\Phi(t)$ is a solution to the matrix equation
$$
\Phi' = [P(t) + Q(t) Y(t)]\Phi, \phh t\in [t_0;T).
$$
By virtue of Liouville formula from the condition $det \Phi(t_0) \ne 0$ of theorem it follows (3.28). Therefore by (2.1) and the uniqueness theorem  $Y(t) = \Psi(t) \Phi^{-1}(t) , \ph t\in[t_0;T)$. Then since $(\Phi(t), \Psi(t))$ is prepared we have $Y(t) = Y^*(t), \ph t\in [t_0;T)$. Hence
$$
y_{12}(t) = y_{21}(t), \phh t\in [t_0;T). \eqno (3.29)
$$
Let $y_1(t)$ and $y_2(t)$ be the solutions to the equations
$$
y' + q_j(t) y^2 + a_{jj}(t) y + F_j(t) = 0, \phh  t\ge t_0, \eqno (3.30)
$$
$$
y' - q_{3-j}(t) y^2 + a_{3-j,3-j}(t) y - F_{3-j}(t) = 0, \phh  t\ge t_0, \eqno (3.31)
$$
respectively with $y_1(t_0) = y_2(t_0) = 0$. By virtue of Theorem 4.1 of work [10] from C), \ph D$_1$) and D$_2$) it follows that $y_1(t)$, \ph $y_2(t)$ exist on $[t_0;+\infty)$ and are nonnegative for all $t\ge t_0$ Moreover if $y_k(t_0) > 0, \ph k=1,2$ then $y_k(t) > 0, \ph t\ge t_0, \ph k=1,2.$. Using Theorem 2.1 to the pairs (3.25), (3.30) and (3.26), (3.32) taking into account (3.27) from here we will get:
$$
y_{11}(t) \ge y_1(t) \ge 0, \phh y_{22}(t) \le - y_2(t) \le 0, \phh t\in [t_0;T), \eqno (3.32)
$$
and if $y_{11}^0 = y_{11}(t_0) >0, \ph y_{22}^0 = y_{22}(t_0) < 0$, then
$$
y_{11}(t) > 0, \phh y_{22}(t) < 0, \phh t\in[t_0;T). \eqno (3.33)
$$
Suppose $T <+\infty$. Then from C) and (3.32) it follows that the function\linebreak $f(t) \equiv \il{T_0}{t} tr [Q(\tau) Y(\tau)]d\tau, \ph t\in[t_0;T)$ is bounded from below on $[t_0;T)$. By Lemma~2.1 from here it follows that $[t_0;T)$ is not the maximum existence interval for $Y(t)$. The obtained contradiction proves (3.24). From (3.24) and (3.28) it follows (3.22), and from (3.24), (3.28) and (3.33) it follows (3.23). The theorem is proved.

\vskip 20pt

\centerline{\bf References}

\vskip 20pt

\noindent
1. L. H. Erbe, Q. Kong and Sh. Ruan, Kamenev type theorems for second order matrix\linebreak \phantom{aa}  differential systems. Proc. Amer. Math. Soc. Vol. 117, Num. 4, 1993, 957 - 962.

\noindent
2. Q. Wang, Oscillation criteria for second order matrix differential systems Arch. Math.\linebreak \phantom{aa} 76 (2001) 385 - 390.

\noindent
3. F. Meng and  A. B. Mingarelli, Oscillation of linear hamiltonian systems, Proc. Amer.\linebreak \phantom{aa} Math. Soc. Vol. 131, Num. 3, 2002, 897 - 904.

\noindent
4. Z. Zhung and S. Zhu, Hartman type oscillation criteria for linear matrix hamiltonian\linebreak \phantom{aa} systems, Dynamic systems and Applications, 17 (2008) 85 -96.

\noindent
5. L. Li, F. Meng and Z. Zhung, Oscillation results related to integral averaging technique\linebreak \phantom{aa} for linear hamiltonian systems, Dynamic systems and Applications 18 (2009) 725 - 736.

\noindent
6. A. B. Mingarelli, On a conjecture for oscillation of second order ordinary differential\linebreak \phantom{aa} systems, Proc. Amer. Math. Soc., Vol. 82. Num. 4, 593 - 598.

\noindent
7. G. J. Butler, L. H. Erbe and A. B. Mingarelly, Riccati techniques and variational\linebreak \phantom{aa} principles in oscillation theory for linear systems, Trans. Amer. Math. Soc. Vol. 303,\linebreak \phantom{aa} Num. 1, 1987, 263 - 231.

\noindent
8. Ph. Hartman, Ordinary differential equations, SIAM - Society for industrial and\linebreak \phantom{aaa} applied Mathematics, Classics in Applied Mathematics 38, Philadelphia 2002.

\noindent
9. G. A. Grigorian.  On two comparison tests for second-order linear  ordinary differential\linebreak \phantom{aa} equations (Russian) Differ. Uravn. 47 (2011), no. 9, 1225 - 1240; translation in Differ.\linebreak \phantom{aa} Equ. 47 (2011), no. 9 1237 - 1252, 34C10.

\noindent
10. G. A. Grigorian, "Two Comparison Criteria for Scalar Riccati Equations with\linebreak \phantom{aa} Applications". Russian Mathematics (Iz. VUZ), 56, No. 11, 17 - 30 (2012).

\noindent
11. G. A. Grigorian, Global Solvability of Scalar Riccati Equations. Izv. Vissh.\linebreak \phantom{aa} Uchebn. Zaved. Mat., 2015, no. 3, pp. 35 - 48.

\noindent
12. G. A. Grigorian, On the Stability of Systems of Two First - Order Linear Ordinary\linebreak \phantom{aa} Differential Equations, Differ. Uravn., 2015, vol. 51, no. 3, pp. 283 - 292.

\noindent
13. G. A. Grigorian. On one oscillatory criterion for the second order linear
 ordinary\linebreak \phantom{aa} differential equations. Opuscula Math. 36, no. 5 (2016), 589–601.
   http://\linebreak \phantom{aaa}\hskip 2pt dx.doi.org/10.7494/OpMath.2016.36.5.589.

\noindent
14. R. Byers, B. J. Harris and M. K. Kwong, Weighted Means and Oscillation Conditions\linebreak \phantom{aa}  for Second Order Matrix Differential Equations. Journal of Differential equations\linebreak \phantom{aa} 61, 164 - 177 (1986).

\end{document}